\documentclass[12pt]{article}
\usepackage{latexsym,amssymb,amsmath}
\begin{document}
\baselineskip=20pt

\newcommand{\la}{\langle}
\newcommand{\ra}{\rangle}
\newcommand{\psp}{\vspace{0.4cm}}
\newcommand{\pse}{\vspace{0.2cm}}
\newcommand{\ptl}{\partial}
\newcommand{\dlt}{\delta}
\newcommand{\sgm}{\sigma}
\newcommand{\al}{\alpha}
\newcommand{\be}{\beta}
\newcommand{\G}{\Gamma}
\newcommand{\gm}{\gamma}
\newcommand{\vs}{\varsigma}
\newcommand{\Lmd}{\Lambda}
\newcommand{\lmd}{\lambda}
\newcommand{\td}{\tilde}
\newcommand{\vf}{\varphi}
\newcommand{\yt}{Y^{\nu}}
\newcommand{\wt}{\mbox{wt}\:}
\newcommand{\rd}{\mbox{Res}}
\newcommand{\ad}{\mbox{ad}}
\newcommand{\stl}{\stackrel}
\newcommand{\ol}{\overline}
\newcommand{\ul}{\underline}
\newcommand{\es}{\epsilon}
\newcommand{\dmd}{\diamond}
\newcommand{\clt}{\clubsuit}
\newcommand{\vt}{\vartheta}
\newcommand{\ves}{\varepsilon}
\newcommand{\dg}{\dagger}
\newcommand{\tr}{\mbox{Tr}}
\newcommand{\ga}{{\cal G}({\cal A})}
\newcommand{\hga}{\hat{\cal G}({\cal A})}
\newcommand{\Edo}{\mbox{End}\:}
\newcommand{\for}{\mbox{for}}
\newcommand{\kn}{\mbox{ker}}
\newcommand{\Dlt}{\Delta}
\newcommand{\rad}{\mbox{Rad}}
\newcommand{\rta}{\rightarrow}
\newcommand{\mbb}{\mathbb}
\newcommand{\lra}{\Longrightarrow}

\begin{center}{\Large \bf Construction of Gel'fand-Dorfman Bialgebras}\end{center}
\begin{center}{\Large \bf from Classical R-Matrices}
\end{center}
\vspace{0.2cm}

\begin{center}{\large Xiaoping Xu}\end{center}
\begin{center}{Department of Mathematics, The Hong Kong University of Science \& Technology}\end{center}
\begin{center}{Clear Water Bay, Kowloon, Hong Kong}\footnote{Research supported
 by Hong Kong RGC Competitive Earmarked Research Grant HKUST6133/00P.}\end{center}
\begin{center}{and Academy of Mathematics \& System Sciences}\end{center}
\begin{center}{Chinese Academy of Sciences, Beijing 10080, P.R.China}\end{center}

\vspace{0.3cm}

\begin{center}{\Large \bf Abstract}\end{center}
\vspace{0.2cm}

{\small Novikov algebras are algebras whose associators are left-symmetric and right multiplication operators are mutually commutative.  A Gel'fand-Dorfman bialgebra is a vector space with a Lie algebra structure and a Novikov algebra structure, satisfying a certain compatibility condition. Such a bialgebraic structure corresponds to a certain Hamiltonian pairs in integrable systems. In this article, we present a construction of Gel'fand-Dorfman bialgebras from certain classical R-matrices on Lie algebras. In particular, we construct such R-matrices from certain abelian subalgebras of Lie algebras. As a result, we show that there exist nontrivial Novikov algebra structures on any finite-dimensional nonzero Lie algebra  over an algebraically closed field with characteristic $0$ or $p>5$ such that they form a Gel'fand-Dorfman bialgebra.} 

\section{Introduction}

A {\it left-symmetric algebra} ${\cal A}$ is a vector space with an algebraic operation $\circ$ satisfying 
$$(u\circ v)\circ w-u\circ (v\circ w)=(v\circ u)\circ w-v\circ (u\circ w)\eqno(1.1)$$ 
for $u,v,w\in{\cal A}$. Left-symmetric algebras play fundamental roles in the theory of affine manifolds (cf. [A], [FD]). A {\it Novikov algebra} is a left-symmetric algebra $({\cal A},\circ)$ satisfying an additional condition
$$(u\circ v)\circ w=(u\circ w)\circ v\qquad\for\;\;u,v,w\in{\cal A}.\eqno(1.2)$$
Novikov and Balinsky [NB] gave the following connection with Poisson brackets of hydrodynamic type
$$\{u(x),v(y)\}=\ptl_x((u\circ v)(x))\dlt(x-y)+(u\circ v+v\circ u)(x)\ptl_x\dlt(x-y).\eqno(1.3)$$

Novikov algebra structure had actually appeared earlier in 
Gel'fand and Dorfman's work [GDo] on Hamiltonian operators. It was named as ``Novikov algebra" by Osborn [O]. A Gel'fand-Dorfman bialgebra is a vector space ${\cal A}$ with two algebraic operations $[\cdot,\cdot]$ and $\circ$ such that $({\cal A},[\cdot,\cdot])$ forms a Lie algebra, $({\cal A},\circ)$ forms a Novikov algebra and the following compatibility condition holds
$$[w,u]\circ v-[w,v]\circ u+[w\circ u,v]-[w\circ v,u]-w\circ [u,v]=0\eqno(1.4)$$
for $u,v,w\in {\cal A}$. Such a bialgebraic structure was introduced by Gel'fand and Dorfman [GDo] in studying certain Hamiltonian pairs, which play important roles in the complete integrability of nonlinear evolution partial differential equations. They correspond to the following Poisson brackets of dynamic type
$$\{u(x),v(y)\}=[u,v](x)\dlt(x-y)+\ptl_x((u\circ v)(x))\dlt(x-y)+(u\circ v+v\circ u)(x)\ptl_x\dlt(x-y).\eqno(1.5)$$
If $({\cal A},[\cdot,\cdot],\circ)$ is a Gel'fand-Dorfman bialgebra, then we say that  $({\cal A},[\cdot,\cdot])$ is a {\it Lie algebra over the Novikov} $({\cal A},\circ)$ and $({\cal A},\circ)$ is a {\it Novikov algebra over the Lie algebra} $({\cal A},[\cdot,\cdot])$. 

The Lie algebras over certain simple Novikov algebras associated with polynomial algebras  were classified by Osborn and Zel'manov [OZ]. In [X3], we gave four general constructions of Gel'fand-Dorfman bialgebras. Moreover, we classified the Lie algebras over more general simple Novikov algebras associated with certain semigroup algebras. In particular, we found that these Lie algebras are generalizations of the well known Block algebras (cf. [B]), which are  simple Lie algebras or ``almost simple'' Lie algebras under certain conditions (cf. [X2]).

For a left-symmetric algebra $({\cal A},\circ)$, we define the  commutator
$$[u,v]=u\circ v-v\circ u\qquad\for\;\;u,v\in{\cal A}.\eqno(1.6)$$
It was stated in [GDo]  that $({\cal A},[\cdot,\cdot])$ forms a Lie algebra and  (1.4) holds (cf. [X3] for the proof). In [X3], we classified all Novikov algebras whose commutator Lie algebras are simple rank-one Witt algebras associated with the group algebras of additive subgroups of the base field with characteristic $0$. 

A linear transformation $T$ on a Lie algebra ${\cal G}$ is called a {\it classical R-matrix} if
$$[T(u),T(v)]=T([T(u),v])+T([u,T(v)])\qquad\for\;\;u,v\in{\cal G}.\eqno(1.7)$$
(cf. [S]). Its relation with the ``classical Yang-Baxter equation'' will be presented in next section for reader's convenience. The main result of this article is as follows.
\psp

{\bf Main Theorem}. {\it Let $T$ be a classical R-matrix on a Lie algebra $({\cal G},[\cdot,\cdot])$ satisfying
$$T[u,T(v)]=[u,T^2(v)]\qquad\mbox{\it for}\;\;u,v\in{\cal G},\eqno(1.8)$$
Then the following algebraic operation
$$u\circ v=T([u,v])-[T(u),v]\qquad\mbox{\it for}\;\;u,v\in{\cal G}\eqno(1.9)$$
defines a Novikov algebra over the Lie algebra ${\cal G}$.}
\psp

In particular, we construct such R-matrices from certain abelian subalgebras of Lie algebras. As a result, we show that there exist nontrivial Novikov algebras over  any finite-dimensional nonzero Lie algebra over an algebraically closed field with characteristic $0$ or $p>5$.

Gel'fand-Dorfman bialgebras are related to certain Lie algebras with one-variable structure, which are natural extensions of the well known loop algebras. 
 Denote by $\mbb{Z}$ the ring of integers and by $\mbb{F}$ a field. All the vector spaces (algebras) in this article are assumed over $\mbb{F}$. Let $({\cal A},[\cdot,\cdot],\circ)$ be a Gel'fand-Dorfman bialgebra. Set
$$\hat{\cal A}={\cal A}\otimes_{\Bbb{F}}\Bbb{F}[t,t^{-1}],\eqno(1.10)$$
where $t$ is an indeterminate. Define an algebraic operation $[\cdot,\cdot]$ on $\hat{\cal A}$ by  
$$[u\otimes t^j,v\otimes t^k]=[u,v]\otimes t^{j+k}+ju\circ v\otimes t^{j+k-1}-kv\circ u\otimes t^{j+k-1}
\eqno(1.11)$$
 for $u,v\in {\cal A},\;j,k\in\Bbb{Z}$. Then $(\hat{\cal A},[\cdot,\cdot])$ forms a Lie algebra. Note that the first term on the right hand side of the above equation defines a loop algebra. The Lie algebra $(\hat{\cal A},[\cdot,\cdot])$ will be useful in the study of vertex algebras, which are fundamental algebraic structures in conformal field theory.

In next section, we shall  prove our main theorem and present a connection with the constant classical Yang-Baxter equation. In Section 3, we shall present a concrete construction based on abelian subalgebras and prove a existence theorem.

\section{General Construction}

In this section, we shall give a general construction of Novikov algebras over any given Lie algebra.
\psp

{\bf Lemma 2.1}. {\it Let ${\cal G}$ be a Lie algebra and let $T$ be a linear transformation on ${\cal G}$. Then the following  algebraic operation $\circ$ on ${\cal G}$:
$$u\circ v=T([u,v])-[T(u),v]\qquad\mbox{\it for}\;\;u,v\in{\cal G},\eqno(2.1)$$
satisfies (1.4).}
\psp

{\it Proof}. For any $u,v,w\in{\cal G}$, we have
\begin{eqnarray*}& &[w,u]\circ v-[w,v]\circ u+[w\circ u,v]-[w\circ v,u]-w\circ [u,v]\\ &=&T([[w,u],v])-[T([w,u]),v]-T([[w,v],u])+[T([w,v]),u]+[T([w, u])\\ & &-[T(w),u],v]-[T([w, v])-[T(w),v],u]-T([w,[u,v]])+[T(w),[u,v]]\\ &=&
T([[w,u],v])-[T([w,u]),v]-T([[w,v],u])+[T([w,v]),u]+[T([w, u]),v]\\ & &-[[T(w),u],v]-[T([w, v]),u]+[[T(w),v],u]-T([w,[u,v]])+[T(w),[u,v]]\\&=&
-[[T(w),u],v]+[[T(w),v],u]+T([[w,u],v])\\ & &-T([[w,v],u])-T([w,[u,v]])+[T(w),[u,v]]\\ &=& [v,[T(w),u]]+[u,[v,T(w)]]+[T(w),[u,v]]\\ & &+T([[w,u],v]+[[v,w],u]+[[u,v],w])=0\hspace{6.6cm}(2.2)\end{eqnarray*}
by skew-symmetry and Jacobi identity of Lie algebra.$\qquad\Box$
\psp

 Let $V$  be a vector space and let
$\{v_j\mid j\in \Omega\}$ be a basis of $V$. Given a linear transformation $R$ on $V\otimes V$, we write
$$R(v_i\otimes v_j)=\sum_{p,q\in\Omega}r_{i,j}^{p,q}v_p\otimes v_q\eqno(2.3)$$
for any $i,j \in \Omega$, where $r_{i,j}^{p,q}\in \mbb{F}$. Define three linear transformations $R^{12},R^{13},R^{23}$ on $V\otimes V\otimes V$ by
$$R^{12}(v_i\otimes v_j\otimes v_k)=\sum_{p,q\in \Omega}r_{i,j}^{p,q}v_p\otimes v_q\otimes v_k,\eqno(2.4)$$
$$R^{13}(v_i\otimes v_j\otimes v_k)=\sum_{p,q\in \Omega}r_{i,k}^{p,q}v_p\otimes v_j\otimes v_q,\eqno(2.5)$$
$$R^{23}(v_i\otimes v_j\otimes v_k)=\sum_{p,q\in \Omega}r_{j,k}^{p,q}v_i\otimes v_p\otimes v_q,\eqno(2.6)$$
for any $i,j,k\in \Omega$. The equation
$$[R^{12},R^{13}]+[R^{13},R^{23}]+[R^{12},R^{23}]=0\eqno(2.7)$$
is called the {\it constant classical Yang-Baxter equation}, where
$[\cdot,\cdot]$ is the commutator of linear transformations.

 Most of the known interesting solutions of the the constant classical Yang-Baxter equation have been constructed by means of Lie algebras. 
Let ${\cal G}$ be a  Lie algebra and let $\{I_j\mid j\in \Omega\}$ be a basis of ${\cal G}$. 
 Suppose that $X$ is an element of  ${\cal G}\otimes {\cal G}$. Write
$$X=\sum_{i,j\in \Omega}r_{i,j}I_i\otimes I_j\eqno(2.8)$$
with $r_{i,j}\in \mbb{F}$. Denote by $U({\cal G})$ the universal
enveloping algebra of ${\cal G}$. Define the elements
$X^{12},\;X^{13},\;X^{23}\in U({\cal G})\otimes U({\cal
G})\otimes U({\cal G})$ by
$$X^{12}=\sum_{i,j\in\Omega}r_{i,j}I_i\otimes I_j\otimes 1,\;
X^{13}=\sum_{i,j\in \Omega}r_{i,j}I_i\otimes 1\otimes I_j,
\;X^{23}=\sum_{i,j\in \Omega}r_{i,j}1\otimes I_i\otimes I_j.\eqno(2.9)$$
Then the {\it constant classical Yang-Baxter equation for} ${\cal G}$ is:
$$[X^{12},X^{13}]+[X^{13},X^{23}]+[X^{12},X^{23}]=0,\eqno(2.10)$$
where $[\cdot,\cdot]$ is the commutator of the tensor algebra $U({\cal G})\otimes U({\cal G})\otimes U({\cal G})$.

Let $V$ be a ${\cal G}$-module with representation $\pi$. If
$X$ is a solution of (2.10), then the linear transformation
$$R=\sum_{i,j\in \Omega}r_{i,j}\pi(I_i)\otimes\pi(I_j)\eqno(2.11)$$
on $V\otimes V$ is a solution of the equation (2.7).

Suppose that the Lie algebra ${\cal G}$ has a nondegenerate invariant symmetric bilinear form $\la\cdot,\cdot\ra$, that is,
$$\la u,v\ra=\la v,u\ra,\qquad \la [u,v],w\ra=\la u,[v,w]\ra\eqno(2.12)$$
for $u,v,w\in{\cal G}$. Take $\{I_j\mid j\in \Omega\}$ to be an orthonormal basis of ${\cal G}$, that is,
$$\la I_i,I_j \ra=\dlt_{i,j}\qquad\for\;\;i,j\in\Omega.\eqno(2.13)$$

For the element $X$ in (2.8), we define a transformation $T_X$ on ${\cal G}$ by
$$T_X(I_i)=\sum_{j\in\Omega}r_{i,j}I_j\qquad\for\;\;i,j\in\Omega.\eqno(2.14)$$
Then the equation (2.10) is equivalent to
$$ [T_X(u),T_X(v)]=T_X([T_X(u),v])+T_X([u,T_X(v)])\qquad\for\;\;u,v\in{\cal G}\eqno(2.15)$$
(e.g. cf. [X1]), which is the condition for $T_X$ to be a classical R-matrix.
\psp

{\bf Theorem 2.2}. {\it Let} $T$ be an R-matrix of a Lie algebra ${\cal G}$. Then the algebraic operation $\circ$ defined in (2.1) gives a Novikov algebra over the Lie algebra ${\cal G}$ if 
$$T([u,v])=[u,T(v)]\qquad\mbox{\it for}\;\;u\in{\cal G},\;v\in T({\cal G}).\eqno(2.16)$$

{\it Proof}. Let $u,v,w\in {\cal G}$. We have
\begin{eqnarray*}& &(u\circ v)\circ w-(u\circ w)\circ v\\ &=&
(T([u, v])-[T(u),v])\circ w-(T([u, w])-[T(u),w])\circ v\\ &=&
T([u, v])\circ w -[T(u),v]\circ w-T([u, w])\circ v+[T(u),w]\circ v\\ &=&
T([T([u, v]), w])-[T^2([u, v]), w]-T([[T(u),v], w])+[T([T(u),v]), w]\\ &&
-T([T([u, w]), v])+[T^2([u, w]), v]+T([[T(u),w], v])-[T([T(u),w]), v]\\&=& 
T([T(u),[w, v]]) +[T([T(u),v]), w]-[T([T(u),w]), v]\\ &=&
[T^2(u),[w, v]]) +[[T^2(u),v], w]-[[T^2(u),w], v]\\ &=&
[T^2(u),[w, v]]) +[w,[v,T^2(u)]]+[v,[T^2(u),w]]=0,\hspace{4.7cm}(2.17)\end{eqnarray*}
by (2.16), and skew symmetry and Jacobi identity of Lie algebra. Note
\begin{eqnarray*} & &u\circ(v\circ w)\\&=&T([u,v\circ w])-[T(u),v\circ w]
\\ &=&T([u,T([v, w])-[T(v),w]])-[T(u),T([v, w])-[T(v),w]]\\ &=& T([u,T([v, w])])-T([u,[T(v),w]])-[T(u),T([v, w])]+[T(u),[T(v),w]].\hspace{1cm}(2.18)\end{eqnarray*}
By (2.15)-(2.18), and skew symmetry and Jacobi identity of Lie algebra, we obtain
\begin{eqnarray*}& &(u\circ v)\circ w-u\circ(v\circ w)-(v\circ u)\circ w+v\circ(u\circ w)\\ &=&T([T([u, v]), w])-[T^2([u, v]), w]-T([[T(u),v], w])+[T([T(u),v]), w]\\ & & -T([u,T([v, w])])+T([u,[T(v),w]])+[T(u),T([v, w])]-[T(u),[T(v),w]]\\ & &-T([T([v, u]), w])+[T^2([v, u]), w]+T([[T(v),u], w])-[T([T(v),u]), w]\\ & & +T([v,T([u, w])])-T([v,[T(u),w]])-[T(v),T([u, w])]+[T(v),[T(u),w]]\\
&=&[T(u),T([v, w])]-T([[T(u),v], w])-T([v,[T(u),w]])-T([u,T([v, w]])\\ & &+T([[T(v),u], w])+T([u,[T(v),w]])+T([v,T([u, w])])-[T(v),T([u, w])]\\ & &
+[T([T(u),v]), w]-[T([T(v),u]), w]+[T(v),[T(u),w]] -[T(u),[T(v),w]]\\&=&
[T(u),T([v, w])]-T([T(u),[v, w]])-T([u,T([v, w]])\\ & &+T([T(v),[u, w]])+T([v,T([u, w])])-[T(v),T([u, w])]\\ & &+[T([T(u),v]), w]-[T([T(v),u]), w]+[[T(v),T(u)],w]=0.\hspace{3.2cm}(2.19)\end{eqnarray*}
Thus $({\cal G},\circ)$ forms a Novikov algebra by (1.1) and (1.2). The above lemma shows that $({\cal G},[\cdot,\cdot],\circ)$ forms a Gel'fand-Dorfman bialgebra.$\qquad\Box$

\section{Existence}

In this section, we shall present examples of Novikov algebras over a given Lie algebras based on abelian subalgebras.
\psp

Let ${\cal G}$ be a nonzero  Lie algebra. Take a subspace $K$ of ${\cal G}$ such that
$$[K,K]=\{0\},\;\;K+[K,{\cal G}]\neq {\cal G}.\eqno(3.1)$$
Pick any  subspace $\bar{K}$ of ${\cal G}$ such that
$${\cal G}=\bar{K}\oplus (K+[K.{\cal G}]).\eqno(3.2)$$
Let $T_0:\bar{K}\rta K$ be any linear map. We define a linear transformation $T$ on ${\cal G}$ by
$$T(u_1+u_2)=T_0(u_1)\qquad \for\;\;u_1\in\bar{K},\;u_2\in K+[K,{\cal G}].\eqno(3.3)$$
Then the linear transformation $T$ trivially satisfies (1.7) and (1.8), Thus (1.9) define a Novikov algebra over ${\cal G}$ by our main theorem. Moreover, we have:
\psp

{\bf Theorem 3.1}. {\it There exist nontrivial Novikov algebras over  any finite-dimensional nonzero Lie algebra over an algebraically closed field with characteristic $0$ or $p>5$.}
\psp

{\it Proof}. Let ${\cal G}$ be any finite-dimensional nonzero Lie algebra.
\pse

{\bf Case 1}. ${\cal G}$ is abelian, that is, $[{\cal G},{\cal G}]=\{0\}$.
\psp

Take any nonzero linear map $\sgm: {\cal G}\rta \mbb{F}$. Define an algebraic operation on ${\cal G}$ by
$$u\circ v=\sgm(v)u\qquad\for\;\;u,v\in{\cal G}.\eqno(3.4)$$
Then for $u,v,w$, we have
$$(u\circ v)\circ w=\sgm(v)u\circ w=\sgm(v)\sgm(w)u,\eqno(3.5)$$
$$u\circ (v\circ w)=\sgm(w)u\circ v=\sgm(w)\sgm(v)u.\eqno(3.6)$$
Thus $({\cal G},\circ)$ is an associative algebra satisfying (1.2). Hence it is a nontrivial Novikov algebra over the abelian Lie algebra  ${\cal G}$.
\pse

{\bf Case 2}. ${\cal G}$ is not abelian.
\pse

First if
$$\mbb{F}u+[u,{\cal G}]={\cal G}\qquad\mbox{for any}\;\;0\neq u\in {\cal G},\eqno(3.7)$$
then ${\cal G}$ is a simple Lie algebra. Classification of finite-dimensional simple Lie algebras over an algebraically closed field with characteristic $0$ or $p>5$ has been done, by which there does not exist a finite-dimensional simple Lie algebra ${\cal G}$ such that (3.7) holds.  If all non-central element $u$ of ${\cal G}$ satisfies (3.7), then ${\cal G}/\mbox{Center}\:{\cal G}$ is a simple Lie algebra satisfies (3.7), which contradicts the classification. Thus there exists a non-central element $u_0\in {\cal G}$ such that 
$$\mbb {F}u_0+[u_0,{\cal G}]\neq {\cal G}.\eqno(3.8)$$

Suppose
$$[u,v]=0\qquad \mbox{for any}\;\; v\in {\cal G}\setminus \mbb{F}u_0+[u_0,{\cal G}].\eqno(3.9)$$
Then there exists $w\in [u_0,{\cal G}]$ such that $[u_0,w]\neq 0$. Take any $u_1\in {\cal G}\setminus (\mbb{F}u_0+[u_0,{\cal G}])$. We have
$$[u_0,u_1+w]=[u_0,u_1]+[u_0,w]=[u_0,w]\neq 0,\;\;u_1+w\not\in  \mbb{F}u_0+[u_0,{\cal G}],\eqno(3.10)$$
which contradicts (3.9). Thus there always exists an element $v_0\in {\cal G}\setminus  (\mbb{F}u_0+[u_0,{\cal G}])$ such that $[u_0,v_0]\neq 0$. 

Take any subspace $V$ of ${\cal G}$ such that
$$\mbb{F}u_0+[u_0,{\cal G}]\subset V,\qquad V\oplus \mbb{F}v_0={\cal G}.\eqno(3.11)$$
We define a linear map $T$ on ${\cal G}$ by
$$T(v_0)=u_0,\;\;T(V)=\{0\}.\eqno(3.12)$$
Then $T$ satisfies (1.7) and (1.8). So (1.9) define a Novikov algebra over ${\cal G}$. Moreover,
$$v_0\circ v_0=T([v_0,v_0])-[T(v_0),v_0]=-[u_0,v_0]\neq 0.\eqno(3.13)$$
Hence (1.9) gives a nontrivial  Novikov algebra over ${\cal G}.\qquad\Box$
\psp

{\bf Remark 3.2}. The real Lie algebra $so(3,\mbb{R})$ satisfies (3.7). Existence of nontrivial Novikov algebras over $so(3,\mbb{R})$ is an interesting open problem.
\psp

{\bf Example}. Let ${\cal G}$ be a finite-dimensional simple Lie algebra over an algebraically closed field $\mbb{F}$ with characteristic 0. Then ${\cal G}$ has a Cartan root decomposition
$${\cal G}=H\oplus\bigoplus_{\al\in\Dlt}\mbb{F}x_\al,\eqno(3.14)$$
where $H$ is a Cartan subalgebra of ${\cal G}$, $\Dlt$ is the set of roots (root system) and each $x_\al$ is  a root vector corresponding to $\al$ (e.g., cf. [H]). Note that any proper subalgebra $K$ of $H$ satisfies (3.1). 

Let $\Dlt_+$ be the subset of positive roots in $\Dlt$. Then any abelian subalgebra $K$ contained in $\sum_{\al\in\Dlt_+}\mbb{F}x_\al$ also satisfies (3.1). 
One can  pick relatively large such abelian subalgebras as follows. It is known that  $\Dlt$ has a subset $\{\al_1,\al_2,...,\al_n\}$, so-called {\it simple roots}, such that any root
$$\al=\sum_{i=1}^n m_i\al_i\;\;\mbox{with}\;\;m_i\in\mbb{Z}\eqno(3.15)$$
Define the {\it height} 
$$\iota(\al)=\sum_{i=1}^nm_i\qquad\mbox{for any root}\;\;\al=\sum_{i=1}^n m_i\al_i.\eqno(3.16)$$
Let
$$k=\mbox{max}\:\{\iota(\al)\mid \al\in\Dlt\}.\eqno(3.17)$$
Set
$$ {\cal G}_m=\sum_{\al\in\Dlt,\;\iota(\al)=m}\mbb{F}x_\al\eqno(3.18)$$
for $0\neq m\in\mbb{Z}$ with $|m|\leq k$. Denote
$$\ell=[\!|k/2|\!]+1.\eqno(3.19)$$
We define
$$K=\sum_{i=\ell}^k{\cal G}_i.\eqno(3.20)$$
Then $K$ is a abelian subalgebra. Moreover,
$$(\sum_{i=-k}^{\ell-k-1}{\cal G}_i)\bigcap (K+[K,{\cal G}])=\{0\}.\eqno(3.21)$$

The following is a simplest example of nontrivial Novikov algebra constructed
from the subalgebra $K$ in (3.20). The three dimensional simple Lie algebra $sl(2,\mbb{F})$ has a basis $\{e_+,e_-,h\}$ with the Lie bracket
$$[h,e_\pm]=\pm 2e_\pm,\qquad [e_+,e_-]=h.\eqno(3.22)$$
We define a linear transformation $T$ on $sl(2,\mbb{F})$ by
$$T(e_-)=e_+,\;\;T(h)=T(e_+)=0.\eqno(3.23)$$
Then (2.1) gives
$$h\circ h=h\circ e_+=e_+\circ h=e_+\circ e_-=e_-\circ e_+=e_+\circ e_+=0,\eqno(3.24)$$
$$h\circ e_-=-2e_+,\;\;e_-\circ h=4e_+,\;\;e_-\circ e_-=-h.\eqno(3.25)$$
The above two expressions give a nontrivial Novikov algebra over $sl(2,\mbb{F})$.
\vspace{1cm}

\noindent{\Large \bf References}

\hspace{0.5cm}

\begin{description}

\item[{[A]}] L. Ausland, Simple transitive groups of affine motions, {\it Amer. J. Math.} {\bf 99} (1977), 809-826.

\item[{[BN]}] A. A. Balinskii and S. P. Novikov, Poisson brackets of hydrodynamic type, Frobenius algebras and Lie algebras, {\it Soviet Math. Dokl.} Vol. {\bf 32} (1985), No. {\bf 1}, 228-231.

\item[{[B]}] R. Block, On torsion-free abelian groups and Lie algebras, {\it Proc. Amer. Math. Soc.} {\bf 9} (1958), 613-620.

\item[{[FG]}] D. Fried and W. Goldman, Three dimensional affine crystallographic groups, {\it Adv. Math}. {\bf 47} (1983), 1-49.

\item[{[GDo]}] I. M. Gel'fand and I. Ya. Dorfman, Hamiltonian operators and algebraic structures related to them, {\it Funkts. Anal. Prilozhen}  {\bf 13} (1979), 13-30.

\item[{[H]}] J. E. Humphreys, {\it Introduction to Lie Algebras and Representation Theory}, Springer-Verlag New York Inc., 1972.

\item[{[O]}]
J. Marshall Osborn, Novikov algebras, {\it Nova J. Algebra} \& {\it Geom.} {\bf 1} (1992), 1-14.

\item[{[OZ]}] J. Marshall Osborn and E. Zel'manov, Nonassociative algebras related to Hamiltonian operators in the formal calculus of variations, {\it J. Pure. Appl. Algebra} {\bf 101} (1995), 335-352.

\item[{[S]}]M. A. Semenov-Tyan-Shanskii,
What is a classical r-matrix? {\it Funct. Anal. Appl.} V. {\bf 17}
(1983), 159-272.

\item[{[X1]}] X. Xu, An analogue of the classical Yang-Baxter equation for general algebraic structures, {\it Mh. Math.} {\bf 119} (1995), 327-346.

\item[{[X2]}] X. Xu, Generalizations of the Block algebras, {\it Manuscripta Math.} {\bf 100} (1999), 489-518.

\item[{[X3]}] X. Xu,  Quadratic conformal superalgebras, {\it J. Algebra} {\bf 231} (2000), 1-38.

\item[{[Z]}]E. I. Zel'manov, On a class of local translation invariant Lie algebras, {\it Soviet Math. Dokl.} Vol {\bf 35} (1987), No. {\bf 1}, 216-218.

\end{description}

\end{document}